\documentclass[11pt]{article}

\begin{document}



\newcommand{\newc}{\newcommand}


\renewcommand{\theequation}{\thesection.\arabic{equation}}
\newc{\eqnoset}{\setcounter{equation}{0}}
\newcommand{\myref}[2]{#1~\ref{#2}}

\newcommand{\mref}[1]{(\ref{#1})}
\newcommand{\reflemm}[1]{Lemma~\ref{#1}}
\newcommand{\refrem}[1]{Remark~\ref{#1}}
\newcommand{\reftheo}[1]{Theorem~\ref{#1}}
\newcommand{\refdef}[1]{Definition~\ref{#1}}
\newcommand{\refcoro}[1]{Corollary~\ref{#1}}
\newcommand{\refprop}[1]{Proposition~\ref{#1}}
\newcommand{\refsec}[1]{Section~\ref{#1}}
\newcommand{\refchap}[1]{Chapter~\ref{#1}}

\newcommand{\beq}{\begin{equation}}
\newcommand{\eeq}{\end{equation}}
\newcommand{\beqno}[1]{\begin{equation}\label{#1}}

\newcommand{\barr}{\begin{array}}
\newcommand{\earr}{\end{array}}

\newc{\bearr}{\begin{eqnarray*}}
\newc{\eearr}{\end{eqnarray*}}

\newc{\bearrno}[1]{\begin{eqnarray}\label{#1}}
\newc{\eearrno}{\end{eqnarray}}

\newc{\non}{\nonumber}
\newc{\nol}{\nonumber\nl}

\newcommand{\bdes}{\begin{description}}
\newcommand{\edes}{\end{description}}
\newc{\benu}{\begin{enumerate}}
\newc{\eenu}{\end{enumerate}}
\newc{\btab}{\begin{tabular}}
\newc{\etab}{\end{tabular}}



\newtheorem{theorem}{Theorem}[section]
\newtheorem{defi}[theorem]{Definition}
\newtheorem{lemma}[theorem]{Lemma}
\newtheorem{rem}[theorem]{Remark}
\newtheorem{exam}[theorem]{Example}
\newtheorem{propo}[theorem]{Proposition}
\newtheorem{corol}[theorem]{Corollary}

\renewcommand{\thelemma}{\thesection.\arabic{lemma}}

\newcommand{\btheo}[1]{\begin{theorem}\label{#1}}
\newc{\brem}[1]{\begin{rem}\label{#1}\em}
\newc{\bexam}[1]{\begin{exam}\label{#1}\em}
\newc{\bdefi}[1]{\begin{defi}\label{#1}}
\newcommand{\blemm}[1]{\begin{lemma}\label{#1}}
\newcommand{\bprop}[1]{\begin{propo}\label{#1}}
\newcommand{\bcoro}[1]{\begin{corol}\label{#1}}
\newcommand{\etheo}{\end{theorem}}
\newcommand{\elemm}{\end{lemma}}
\newcommand{\eprop}{\end{propo}}
\newcommand{\ecoro}{\end{corol}}
\newc{\erem}{\end{rem}}
\newc{\eexam}{\end{exam}}
\newc{\edefi}{\end{defi}}

\newc{\rmk}[1]{{\bf REMARK #1: }}
\newc{\DN}[1]{{\bf DEFINITION #1: }}

\newcommand{\bproof}{{\bf Proof:~~}}
\newc{\eproof}{{\vrule height8pt width5pt depth0pt}\vspace{3mm}}


\newcommand{\rarr}{\rightarrow}
\newcommand{\Rarr}{\Rightarrow}
\newcommand{\tru}{\backslash}
\newc{\bfrac}[2]{\dspl{\frac{#1}{#2}}}


\newc{\nl}{\vspace{2mm}\\}
\newc{\nid}{\noindent}


\newcommand{\oneon}[1]{\frac{1}{#1}}
\newcommand{\dspl}{\displaystyle}
\newc{\grad}{\nabla}
\newc{\Div}{\mbox{div}}
\newc{\pdt}[1]{\dspl{\frac{\partial{#1}}{\partial t}}}
\newc{\pdn}[1]{\dspl{\frac{\partial{#1}}{\partial \nu}}}
\newc{\pdNi}[1]{\dspl{\frac{\partial{#1}}{\partial \mathcal{N}_i}}}
\newc{\pD}[2]{\dspl{\frac{\partial{#1}}{\partial #2}}}
\newc{\dt}{\dspl{\frac{d}{dt}}}
\newc{\bdry}[1]{\mbox{$\partial #1$}}
\newc{\sgn}{\mbox{sign}}

\newc{\Hess}[1]{\frac{\partial^2 #1}{\pdh z_i \pdh z_j}}
\newc{\hess}[1]{\partial^2 #1/\pdh z_i \pdh z_j}


\newcommand{\Coone}[1]{\mbox{$C^{1}_{0}(#1)$}}
\newcommand{\lspac}[2]{\mbox{$L^{#1}(#2)$}}
\newc{\hspac}[2]{\mbox{$C^{0,#1}(#2)$}}
\newc{\Hspac}[2]{\mbox{$C^{1,#1}(#2)$}}
\newc{\Hosp}{\mbox{$H^{1}_{0}$}}
\newcommand{\Lsp}[1]{\mbox{$L^{#1}(\Og)$}}
\newc{\hsp}{\Hosp(\Og)}


\newc{\ag}{\alpha}
\newc{\bg}{\beta}
\newc{\cg}{\gamma}\newc{\Cg}{\Gamma}
\newc{\dg}{\delta}\newc{\Dg}{\Delta}
\newc{\eg}{\varepsilon}
\newc{\zg}{\zeta}
\newc{\thg}{\theta}
\newc{\llg}{\lambda}\newc{\LLg}{\Lambda}
\newc{\kg}{\kappa}
\newc{\rg}{\rho}
\newc{\sg}{\sigma}\newc{\Sg}{\Sigma}
\newc{\tg}{\tau}
\newc{\fg}{\phi}\newc{\Fg}{\Phi}
\newc{\vfg}{\varphi}
\newc{\og}{\omega}\newc{\Og}{\Omega}
\newc{\pdh}{\partial}

\newc{\ccG}{{\cal G}}


\newc{\ii}[1]{\int_{#1}}
\newc{\iidx}[2]{{\dspl\int_{#1}~#2~dx}}
\newc{\bii}[1]{{\dspl \ii{#1} }}
\newc{\biii}[2]{{\dspl \iii{#1}{#2} }}
\newc{\su}[2]{\sum_{#1}^{#2}}
\newc{\bsu}[2]{{\dspl \su{#1}{#2} }}

\newcommand{\iiomdx}[1]{{\dspl\int_{\Og}~ #1 ~dx}}
\newc{\biiom}[1]{{\dspl\int_{\bdrom}~ #1 ~d\sg}}
\newc{\io}[1]{{\dspl\int_{\Og}~ #1 ~dx}}
\newc{\bio}[1]{{\dspl\int_{\bdrom}~ #1 ~d\sg}}
\newc{\bsir}{\bsu{i=1}{r}}
\newc{\bsim}{\bsu{i=1}{m}}

\newc{\iibr}[2]{\iidx{\bprw{#1}}{#2}}
\newc{\Intbr}[1]{\iibr{R}{#1}}
\newc{\intbr}[1]{\iibr{\rg}{#1}}
\newc{\intt}[3]{\int_{#1}^{#2}\int_\Og~#3~dxdt}

\newc{\itQ}[2]{\dspl{\int\hspace{-2.5mm}\int_{#1}~#2~dz}}
\newc{\mitQ}[2]{\dspl{\rule[1mm]{4mm}{.3mm}\hspace{-5.3mm}\int\hspace{-2.5mm}\int_{#1}~#2~dz}}
\newc{\mitQQ}[3]{\dspl{\rule[1mm]{4mm}{.3mm}\hspace{-5.3mm}\int\hspace{-2.5mm}\int_{#1}~#2~#3}}

\newc{\mitx}[2]{\dspl{\rule[1mm]{3mm}{.3mm}\hspace{-4mm}\int_{#1}~#2~dx}}

\newc{\mitQq}[2]{\dspl{\rule[1mm]{4mm}{.3mm}\hspace{-5.3mm}\int\hspace{-2.5mm}\int_{#1}~#2~d\bar{z}}}
\newc{\itQq}[2]{\dspl{\int\hspace{-2.5mm}\int_{#1}~#2~d\bar{z}}}

\newc{\pder}[2]{\dspl{\frac{\partial #1}{\partial #2}}}


\newc{\ui}{u_{i}}
\newcommand{\upl}{u^{+}}
\newcommand{\umn}{u^{-}}
\newcommand{\un}{\{ u_{n}\}}

\newcommand{\uo}{u_{0}}
\newc{\voi}{v_{i}^{0}}
\newc{\uoi}{u_{i}^{0}}
\newc{\vu}{\vec{u}}

\newc{\xo}{x_{0}}
\newc{\Br}{B_{R}}
\newc{\Bro}{\Br (\xo)}
\newc{\bdrom}{\bdry{\Og}}
\newc{\ogr}[1]{\Og_{#1}}
\newc{\Bxo}{B_{x_0}}

\newc{\inP}[2]{\|#1(\bullet,t)\|_#2\in\cP}
\newc{\cO}{{\mathcal O}}
\newc{\inO}[2]{\|#1(\bullet,t)\|_#2\in\cO}

\newc{\newl}{\\ &&}

\newc{\bilhom}{\mbox{Bil}(\mbox{Hom}(\RR^{nm},\RR^{nm}))}
\newc{\VV}[1]{{V(Q_{#1})}}

\newc{\ccA}{{\mathcal A}}
\newc{\ccB}{{\mathcal B}}
\newc{\ccC}{{\mathcal C}}
\newc{\ccD}{{\mathcal D}}
\newc{\ccE}{{\mathcal E}}
\newc{\ccH}{\mathcal{H}}
\newc{\ccF}{\mathcal{F}}
\newc{\ccI}{{\mathcal I}}
\newc{\ccJ}{{\mathcal J}}
\newc{\ccP}{{\mathcal P}}
\newc{\ccQ}{{\mathcal Q}}
\newc{\ccR}{{\mathcal R}}
\newc{\ccS}{{\mathcal S}}
\newc{\ccT}{{\mathcal T}}
\newc{\ccX}{{\mathcal X}}
\newc{\ccY}{{\mathcal Y}}
\newc{\ccZ}{{\mathcal Z}}

\newc{\bb}[1]{{\mathbf #1}}
\newc{\bbA}{{\mathbf A}}
\newc{\myprod}[1]{\langle #1 \rangle}
\newc{\mypar}[1]{\left( #1 \right)}


\newc{\lspn}[2]{\mbox{$\| #1\|_{\Lsp{#2}}$}}
\newc{\Lpn}[2]{\mbox{$\| #1\|_{#2}$}}
\newc{\Hn}[1]{\mbox{$\| #1\|_{H^1(\Og)}$}}


\newc{\cyl}[1]{\og\times \{#1\}}
\newc{\cyll}{\og\times[0,1]}
\newc{\vx}[1]{v\cdot #1}
\newc{\vtx}[1]{v(t,x)\cdot #1}
\newc{\vn}{\vx{n}}

\newcommand{\RR}{{\rm I\kern -1.6pt{\rm R}}}


\newenvironment{proof}{\noindent\textbf{Proof.}\ }
{\nopagebreak\hbox{ }\hfill$\Box$\bigskip}


\newc{\itQQ}[2]{\dspl{\int_{#1}#2\,dz}}
\newc{\mmitQQ}[2]{\dspl{\rule[1mm]{4mm}{.3mm}\hspace{-4.3mm}\int_{#1}~#2~dz}}
\newc{\MmitQQ}[2]{\dspl{\rule[1mm]{4mm}{.3mm}\hspace{-4.3mm}\int_{#1}~#2~d\mu}}

\newc{\MUmitQQ}[3]{\dspl{\rule[1mm]{4mm}{.3mm}\hspace{-4.3mm}\int_{#1}~#2~d#3}}
\newc{\MUitQQ}[3]{\dspl{\int_{#1}~#2~d#3}}

\vspace*{-.8in}
\begin{center} {\LARGE\em Global Existence and Global Attractors  of Cross Diffusion Systems on Planar Domains.}

 \end{center}

\vspace{.1in}

\begin{center}

{\sc Dung Le}{\footnote {Department of Mathematics, University of
Texas at San
Antonio, One UTSA Circle, San Antonio, TX 78249. {\tt Email: Dung.Le@utsa.edu}\\
{\em
Mathematics Subject Classifications:} 35J70, 35B65, 42B37.
\hfil\break\indent {\em Key words:} Cross diffusion systems,  H\"older
regularity, global existence.}}

\end{center}

\begin{abstract}
Global existence of strong solutions and the existence of global and atrractors are established for generalized Shigesada-Kawasaki-Teramoto models on planar domains. The cross diffusion and reaction can have polynomial growth of any order.
\end{abstract}

\vspace{.2in}

\newc{\BLLg}{\mathbf{\LLg}}

\section{Introduction}\label{introsec}\eqnoset
Shigesada {\it et al.} in \cite{SKT} introduce the following model
\beqno{e0}\left\{\barr{lll} u_t &=& \Div[\nabla(a_1u+\ag_{11}u^2+\ag_{12}uv)+b_1u\nabla\Fg(x)]+f_1(u,v),\\v_t &=& \Div[\nabla(a_2v+\ag_{21}uv+\ag_{22}v^2)+b_2v\nabla\Fg(x)]+f_2(u,v),\earr\right.\eeq where $f_i(u,v)$ are reaction terms of Lotka-Volterra type and quadratic in $u,v$. The unknowns $u(x,t),v(x,t)$ denote the densities of two species at time $t$ and location $x\in\Og$, a bounded domain in $\RR^2$. Dirichlet or Neumann boundary conditions were usually assumed for \mref{e0}. This model was used to describe the population dynamics of the species $u,v$ which move under the influence of population pressures and of the environmental potential $\Fg(x)$.

Under suitable assumptions of the coefficients in \mref{e0}, Yagi proved in \cite{yag} the global existence of solutions to the above system for a planar domain $\Og$. Clearly, \mref{e0} is a special case of the following system
\beqno{e00z}u_t=\Delta(P(u))+\hat{f}(u,Du)\quad  (x,t)\in Q=\Og\times(0,T),\eeq
where $m\ge2$, $u:\Og\to\RR^m$, $P:\RR^m\to\RR^m$, whose components are {\em polynomials} in $u$, and $f:\RR^m\times\RR^{nm}\to\RR^m$ are vector valued functions. The potential $\Fg(x)$ is incorporated in $\hat{f}(u,Du)$, which will be assumed to have {\em linear} growth in $Du$.  In this paper, we then refer to the above system as the {\em generalized (SKT) system} and allow $P,\hat{f}$ to have {\em polynomial growth of any order} in $u$.

Under suitable assumptions of the parameters $\ag_{ij}$'s in \mref{e0}, Yagi proved in \cite{yag} that the solutions with positive initial data will stay positive and the Jacobian of $P(u,v)$ 
$$A(u,v) = \left[\barr{cc} a_1+2\ag_{11}u+\ag_{12}v & \ag_{12}u\\ \ag_{21}v& a_2+\ag_{21}u+2\ag_{22}v\earr\right]$$
is uniformly elliptic for $u,v\ge0$. In fact, there are positive constants $C$ and $c_i$'s, depending on the parameters $d_i$'s and $\ag_{ij}$'s, and a $C^1$ function $\llg(u,v)\sim c_0+c_1u+c_2v$ such that for any $\zeta\in\RR^{4}$ and nonegative $u,v$  we have
$$\llg(u,v)|\zeta|^2 \le \myprod{A(u,v)\zeta,\zeta} \mbox{ and } |A(u,v)|\le C\llg(u,v).$$ 

In this paper, we consider the following generalized  version of \mref{e00z} of $m$ equations ($m\ge2$)
 and rewrite it in a much more general form as 
\beqno{e1}\left\{\barr{ll} u_t=\Div(A(u)Du)+\hat{f}(u,Du)& (x,t)\in Q=\Og\times(0,T),\\u(x,0)=U_0(x)& x\in\Og\\\mbox{$u=0$ or $\frac{\partial u}{\partial \nu}=0$ on $\partial \Og\times(0,T)$}. &\earr\right.\eeq

We will always assume that the initial data $U_0$ is given in $W^{1,p_0}(\Og,\RR^m)$ for some $p_0>2$, the dimension of $\Og$. As usual, $W^{1,p}(\Og,\RR^m)$, $p\ge1$, will denote the standard Sobolev spaces whose elements are vector valued functions $u\,:\,\Og\to \RR^m$ with finite norm $$\|u\|_{W^{1,p}(\Og,\RR^m)} = \|u\|_{L^p(\Og)} + \|Du\|_{L^p(\Og)}.$$

Inspired by \mref{e0} and the above discussion, for $m\ge2$ we assume the following more general condition on the ellipticity of $A(u)$.

\bdes

\item[A)] $A(u)$ is $C^1$ in $u$ and there are constants $\llg_0,C>0$  and a scalar $C^1$ function $\llg(u)$  such that  for all $u\in\RR^m$ and $\zeta\in\RR^{2m}$ 
\beqno{A1} \llg(u)\ge \llg_0,\;\llg(u)|\zeta|^2 \le \myprod{A(u)\zeta,\zeta} \mbox{ and } |A(u)|\le C\llg(u).\eeq 

In addition, $|A_u|\le C|\llg_u|$ and  the following number  is finite:\beqno{LLgx} \mathbf{\LLg}=\sup_{W\in\RR^m}\frac{|\llg_W(W)|}{\llg(W)}.\eeq

\edes 

Here and throughout this paper, if $B$ is a $C^1$ function in $u\in \RR^m$ then we abbreviate it  derivative $\frac{\partial B}{\partial u}$ by $B_u$.

Concerning the reaction term $\hat{f}$, we will assume the following.

\bdes \item[F)] There exist a constant $C$ and a function $f(u)$ such that \beqno{FUDU11}|\hat{f}(u,Du)| \le C\llg^\frac12(u)|Du| + f(u),\eeq \beqno{fUDU21}|f_u(u)| \le C\llg(u).\eeq
\edes

Supposing that $\llg(u)$ {\em is bounded from above}, under the assumption A) and F), the global existence of \mref{e1} was studied in \cite{Am1} for bounded domains $\Og\in\RR^n$, $n\ge2$. It was shown in \cite{Am1} that a solution $u$ of \mref{e1} exists globally if its $W^{1,p}(\Og)$ norm for some $p>n$ does not blow up in finite time. First of all, the assumption that $\llg(u)$ is a constant or bounded does not apply to \mref{e1} in general because maximum principles are not available to show that $u$ is bounded. Even if one knows that $u$ is bounded, only its partial regularity properties is established, see \cite{GiaS}.

In our recent work \cite{dleANS},  estimates for the $W^{1,p}(\Og)$ norms for some $p>n$ of a strong solution and then its global existence were established for \mref{e1} under A), F) and the weakest assumption that this solution is uniformly VMO (Vanishing Mean Oscillation) in the assumption M') of \cite{dleANS}. No boundedness assumption is needed. The proof in \cite{dleANS} relies on fixed point theories, instead of the semigroup approach in \cite{Am1}, and weighted Gagliardo-Ninrehnberg inequalities involving BMO norms. However, the checking of the uniform VMO assumption in \cite{dleANS} is not a simple task.

In this paper, for planar domains $\Og\subset\RR^2$ we will show in \reftheo{extcoro1} that it is sufficient to control the $W^{1,2}(\Og)$ of a strong solution of \mref{e1} to establish its global existence. If this can be done for all initial data $U_0$ in $X=W^{1,p_0}(\Og)$ then \mref{e1} defines a global semiroup $\{\ccS(t)\}_{t\ge0}$ on $X$, namely $$\ccS(0)U_0=U_0,\quad \ccS(t)U_0(x)=u(x,t)$$ is defined for all $t>0$, with $u$ being the solution of \mref{e1}. We will show further in \reftheo{extcoro2} that if the $W^{1,2}(\Og)$ norms of the strong solutions are uniformly bounded for $t$ large then this semigroup possesses a global attractor and exponential attractors in $X$. Let us recall the definition of a global attractor: A set $\ccA\subset X$ is a global (or universal) attractor if 1) $\ccA$ is an invariant set ($S(t)\ccA=\ccA$, $\forall t\ge0$), 2) For any $U_0\in X$ $S(t)U_0$ converges to $\ccA$ as $t\to\infty$.  in the Banach space $X$. The notion of exponential attractors in Hilbert spaces was introduced in \cite{ET} and the same definition applies for Banach spaces.

We state our main results in \refsec{res} and present their proof in \refsec{w12est}. In \refsec{2dsktsec}, we apply our main theorems to the generalized SKT system \mref{e00z}. We will assume that $\llg(u)$ has polynomial growth in $u$, i.e. $\llg(u)\sim \llg_0+\llg_1|u|^k$ for some $k>0$, and the results in our main theorems continues to hold under much weaker assumptions. Namely, one needs only control the $\|u\|_{L^q(\Og)}$ for some $q>k$. In fact, for $k\in(0,2]$ it is sufficient to estimate $\|u\|_{L^1(\Og)}$ (the case $k>2$ needs a mild assumption on the number $C_*$ in A)). We conclude the paper by showing that this is the case if the reaction term of \mref{e00z} is of competitive type in some sense.

\section{Preliminaries and Main Results}\eqnoset\label{res}

We state the main results of this paper in this this section. Their proof will be given in the next section.

Our first result concerns the global existence of strong solutions to \mref{e1}. We imbed \mref{e1} in the following family of systems 
 \beqno{e1famz}\left\{\barr{l} u_t=\Div(A(\sg u)Du)+\hat{f}(\sg u,\sg Du)\quad (x,t)\in Q=\Og\times(0,T_0), \sg\in[0,1]\\u(x,0)=U_0(x)\quad  x\in\Og\\\mbox{$u=0$ or $\frac{\partial u}{\partial \nu}=0$ on $\partial \Og\times(0,T_0)$}. \earr\right.\eeq
 
In applications, $\llg(u)$ usually has polynomial growth in $u$ so that we will introduce the following stronger version of \mref{LLgx}.
\beqno{llgullgeg} |\llg_u(u)|\le \BLLg_1\llg^{1-\eg_0}(u)\quad \mbox{for some $\BLLg_1,\eg_0>0$ and all $u\in\RR^m$}.\eeq

\btheo{extcoro1}  Let $\Og\subset\RR^2$. Suppose A), F), \mref{llgullgeg} and  that we can establish the following: For any $T_0>0$ there is a constant $M(U_0,T_0)$ depending on $\|U_0\|_{W^{1,p_0}(\Og)}$ and $T_0$ such that  any strong solution $u_\sg$ of \mref{e1famz} on $\Og\times(0,T_0)$
\beqno{Du2boundz} \sup_{t\in(0,T_0)}\|u_\sg(\cdot,t)\|_{W^{1,2}(\Og)}\le M(U_0,T_0),\eeq and for some $s_0>0$  \beqno{llgr} \sup_{t\in(0,T_0)}\|\llg^{s_0}(u_\sg(\cdot,t))\|_{L^1(\Og)}\le M(U_0,T_0).\eeq 

Then \mref{e1} has a unique solution which  exists globally.\etheo

\brem{fsigmarem} The assumptions on the bounds for solutions of \mref{e1famz} in \reftheo{extcoro2} are not very restrictive as they seem at first glance. In applications, since the systems in \mref{e1famz} usually satisfy A) and F) for the same set of constants so that we need to check the estimates \mref{Du2boundz} and \mref{llgr} only for $\sg=1$. Furthermore, the reaction terms $\hat{f}(\sg u,\sg Du)$ can also be replaced by $\hat{f}_\sg(u,Du)$ if these functions satisfy F) uniformly.
\erem

Thus, if \mref{Du2boundz} holds for any initial data $U_0 \in X:=W^{1,p_0}(\Og)$ then \mref{e1} defines a global semiroup $\{\ccS(t)\}_{t\ge0}$ on $X$, namely $$\ccS(0)U_0=U_0,\quad \ccS(t)U_0(x)=u(x,t)$$ is defined for all $t>0$, with $u$ being the solution of \mref{e1}.

Next, we will show that if the norm $\|u(\cdot,t)\|_{W^{1,2}(\Og)}$ can be bounded uniformly for bounded initial data $U_0 \in X$  when $t$ is sufficiently large then the global dynamical systems $\{\ccS(t)\}_{t\ge0}$ possesses a global attractor in $X$.

\btheo{extcoro2}  Assume as in \reftheo{extcoro1}.  Suppose further that the dynamical system defined by \mref{e1} possesses an absorbing ball in $W^{1,2}(\Og)$. That is there is a constant $M$ such that for any bounded set $K\subset W^{1,p_0}(\Og)$ there is  $T_K>0$ such that any global solution $u$ of \mref{e1} with initial data $U_0\in K$ will satisfy  \beqno{Du2boundzz}\|u(\cdot,t)\|_{W^{1,2}(\Og)}\le M \quad \mbox{for all $t\ge T_K$},\eeq and for some $s_0>0$ \beqno{llgrzz} \|\llg^{s_0}(u(\cdot,t))\|_{L^1(\Og)}<M \quad \mbox{for all $t\ge T_K$}.\eeq 

Then the system \mref{e1} possesses a global attractor in $X=W^{1,p_0}(\Og)$.\etheo

We remark that the conditions in the above theorem also show that $\{\ccS(t)\}_{t\ge0}$ possesses exponential attractors in the Banach space $X$. The notion of exponential attractors in Hilbert spaces was introduced in \cite{ET}: A set $\ccA\subset X$ is an exponential attractor if 1) $\ccA$ is an positively invariant set ($S(t)\ccA\subset\ccA$, $\forall t\ge0$), 2) For any $U_0\in X$ $S(t)U_0$ converges expontiallly to $\ccA$ as $t\to\infty$. The same definition applies when $X$ is a Banach space. It is shown in \cite{DN} that this notion is quite universal: if $\{\ccS(t)\}_{t\ge0}$ possesses a global attractor and $\ccS(t)$ is a $C^1$ compact map on $X$ for all $t>0$ then exponential attractors exist. In the next section, where we present the proof of the above theorems, higher regularity of $u$ will be established and we will see that $Du$ is H\"older continuous in $(x,t)$ (see \mref{uniDuzzz}). From this, we can prove that the $\ccS(t)$ is $C^1$ on $X$ and thus, by \cite{DN},  $\{\ccS(t)\}_{t\ge0}$ possesses exponential attractors in the Banach space $X$.

\section{Global existence results in the general case and the proof of the main theorems}\eqnoset\label{w12est}

The proof of our main theorems relies on the global existence result in our recent work \cite{dleANS}, which deals with the general case $n\ge2$. In order to describe the main result in \cite{dleANS}, we first recall the following main technical result and introduce the key condition M').

For the general case $n\ge2$, we also assumed in \cite{dleANS} that
 \bdes\item[SG)] (The spectral gap condition)  $(n-2)/n < C_*^{-1}$. \edes

\blemm{dleANSprop}  We assume that $A,\hat{f}$ satisfy A), F) and SG).
Suppose that $u$ is strong solution  of \mref{e1} on $\Og\times(0,T_0)$ and there are a constant $C_{1,2}$, which may depend on $T_0$, and a sufficiently large $r$ such that
 \beqno{Du2bound} \int_0^{T_0}\iidx{\Og}{|Du(x,s)|^2}ds<C_{1,2},\eeq
 \beqno{llgrzzz} \|\llg(u(\cdot,t))\|_{L^r(\Og)}<C_{1,2}.\eeq 
More importantly, we assume that
\bdes\item[M')] 
for any given $\mu_0>0$ there is a positive $R_{\mu_0}$ such that   \beqno{KeyVMO} \mathbf{\LLg}^2\sup_{x_0\in\bar{\Og},t\in(0,T_0)}\|u(\cdot,t)\|_{BMO(B_{R_{\mu_0}}(x_0)\cap\Og)}^2 \le \mu_0.\eeq  
\edes

Then if $\mu_0$ is sufficiently small in terms of the parameters in A) then there are a number $p>n$ and a constant $C_{1,p}$ depending on the parameters in A), F) and $\mathbf{\LLg},\mu_0,R_{\mu_0}$  and the geometry of $\Og$ such that
\beqno{uniDu}\|u(\cdot,t)\|_{W^{1,p}(\Og)} \le C_{1,p}\quad \forall t\in(0,T_0).\eeq
\elemm

This result is a consequence of  \cite[Proposition 3.1]{dleANS}  by taking $\bg(u)=\llg^{-1}(u)$ and $W=U=u$. The conditions required by this proposition that $\|\bg(u(\cdot,t))\|_{L^r(\Og)}$ is bounded and that $\llg(u)\bg(u)$ is a $A_\frac43$ weight are satisfied here because $\llg(u)$ is bounded from below and $\llg(u)\bg(u)=1$.

\brem{Oggeo}The dependence of $C_{1,p}$ in \mref{uniDu} on the geometry of $\Og$ means: $C_{1,p}$ depends on a number $N_{\mu_0}$  of balls $B_{R_{\mu_0}}(x_i)$, $x_i\in\bar{\Og}$ and $R_{\mu_0}$ is as in M'), such that
\beqno{Ogcover} \bar{\Og} \subset \cup_{i=1}^{N_{\mu_0}}B_{R_{\mu_0}}(x_i).\eeq The bound \mref{uniDu} was established locally for balls $B_{R_{\mu_0}}(x_i)\subset\Og$ with $R_{\mu_0}$ satisfying M'). If $x_i$ is on the boundary $\partial\Og$ then a flatenning and odd/even reflection arguments, depending on the type of boundary condition of $u$, can apply to extend the proof for the interior case to the boundary one. Adding these estimates, we obtain the global \mref{uniDu}.  Hence, $N_{\mu_0}$ and $C_{1,p}$ also depend on the geometry of $\partial\Og$. Also, from the proof one can see that $C_{1,p}$ also depends on $R_{\mu_0}^{-1}$.
\erem

\reflemm{dleANSprop} is the key ingredient in the proof of the following theorem, which is a consequence of  \cite[Theorem 4.2]{dleANS}  by taking $\bg(u)=\llg^{-1}(u)$ again. 

\btheo{extthm} Suppose that \mref{e1} satisfies A), F). We consider the following family of systems 
\beqno{e1fam}\left\{\barr{l} u_t=\Div(A(\sg u)Du)+\hat{f}(\sg u,\sg Du)\quad (x,t)\in Q=\Og\times(0,T_0), \sg\in[0,1]\\u(x,0)=U_0(x)\quad  x\in\Og\\\mbox{$u=0$ or $\frac{\partial u}{\partial \nu}=0$ on $\partial \Og\times(0,T_0)$}. \earr\right.\eeq
 Assume that the strong solutions of these systems satisfy \mref{Du2bound}, \mref{llgrzzz} and M') uniformly. 
Then \mref{e1} has a unique strong solution $u$ which exists globally on $\Og\times(0,\infty)$.

\etheo

There is a subtle point in the assumptions in \cite[Theorem 4.2]{dleANS} that is worth discussing here. The paper \cite{dleANS} dealt with the general case where $n\ge2$, $A,\hat{f}$ depend also on $x,t$ and so does $\llg$. It then assumed that $\llg(x,t,u)$ is bounded near $t=0$. It turns out that this is not needed for the global existence result. In this paper we assume that $A,\hat{f}$ are independent of $x,t$ so that we need to explain this matter further by sketching the main ideas in proof of \cite[Theorem 4.2]{dleANS} below.

The proof of \cite[Theorem 4.2]{dleANS}, or its special case \reftheo{extthm} here, makes use of fixed point theorems in a very standard way by considering the linear compact maps associated to the systems \mref{e1fam} on the Banach space $\ccX=\ccX_1\cap\ccX_3$ (the space $\ccX_2=C(\Og\times(0,t_0))$ in the proof of \cite[Theorem 4.2]{dleANS} is not needed here) where
 $$\ccX_1=C((0,T_0),C^{\ag_0}(\Og))\mbox{ and } \ccX_3=\{u:\, Du\in C^{\ag_0,\ag_0/2}(\Og\times(t_0,T_0))\}.$$ Here, $t_0\in(0,T_0)$ is fixed and  $\ag_0>0$ is a number such that $W^{1,p}(\Og)\cap W^{1,p_0}(\Og)$, with $p>n$ being given in \mref{uniDu}, is compactly imbedded in $C^{\ag_0}(\Og))$. We define $\|u\|_{\ccX}=\|u\|_{\ccX_1}+\|u\|_{\ccX_3}$, where
 $$\|u\|_{\ccX_1}=\sup_{t\in(0,T_0)}\|u\|_{C^{\ag_0}(\Og)}\mbox{ and }\|u\|_{\ccX_3}=\sup_{t\in(t_0,T_0)}\|Du\|_{C^{\ag_0,\ag_0/2}(\Og)}.$$

 For each $w\in\ccX$ and $\sg\in[0,1]$ we define $u=T_\sg(w)$ to be the weak solution of \beqno{e1famlin}\left\{\barr{l} u_t=\Div(A(\sg w)Du)+\hat{f}(\sg w,\sg Dw)\quad (x,t)\in Q=\Og\times(0,T_0), \sg\in[0,1]\\u(x,0)=U_0(x)\quad  x\in\Og\\\mbox{$u=0$ or $\frac{\partial u}{\partial \nu}=0$ on $\partial \Og\times(0,T_0)$}. \earr\right.\eeq
 
  The strong solution of \mref{e1} on the cylinder $Q$ is then the fixed point of the corresponding map, $\sg=1$, in $\ccX$.

 From the regularity theory of linear systems with smooth data (see \cite[Chapter 4]{lieb} or \cite{AF}), we see that the solution of the above system is in $\ccX$. Furthermore, the higher regularity estimates in \cite{AF,lieb} also show that $T_\sg$ is a compact map on $\ccX$.

 We now discuss the uniform boundedness of the fixed points $u=T_\sg(u)$. First of all, we observe that $u$ is a strong solution on $(0,T_0)$. In fact, as $u\in\ccX$ we have that $Du$ is H\"older continuous in $Q$. Using the regularity theory of linear systems with smooth data again, we see that $u$ is a strong solution in $\Og\times I$. Therefore, \reflemm{dleANSprop} is applicable here.

The argument in \cite[Proposition 3.1]{dleANS}, which gives \reflemm{dleANSprop},  made use of a cutoff function $\eta$ for the interval $[0,T+t_0]$ and $[0,T+2t_0]$, that is $\eta(t)=0$ for $t<T+t_0$ and $\eta(t)=1$ for $t>T+2t_0$, to avoid the dependence on the initial data at $t=0$. If we allow $t_0=0$ then the bound for $\|u(\cdot,t)\|_{W^{1,p}(\Og)}$ in \cite[Proposition 3.1]{dleANS} (see also \cite[inequality 3.13]{dleANS}) is independent of $t_0$ but $\|u(\cdot,0)\|_{W^{1,p}(\Og)}$. Of course, we can take $p\in (n,p_0]$ so that \mref{uniDu} holds for all $t\in[0,T_0)$ with $C_{1,p}$ explicitly depending on $\|U_0\|_{W^{1,p_0}(\Og)}$. We then obtain a uniform bound for $\|u(\cdot,t)\|_{W^{1,p}(\Og)}$, $t\in(0,T_0)$. Hence, $\|u\|_{\ccX_1}$ is uniformly bounded for any fixed points $u$ of $T_\sg$.

Now, for any strong solution $u$ to \mref{e1fam} and cylinder $Q_R=B_R\times(t-R^2,t)$ in $Q$ we can use H\"older's inequality (in the $x$ integral) and \mref{uniDu} to obtain ($dz=dxdt$)
\beqno{giascond} R^{-n}\itQ{Q_R}{|Du|^2}\le R^{-n+2+n(1-\frac2p)}C_{1,p}^\frac2p=R^{2(1-\frac{n}{p})}C_{1,p}^\frac2p.\eeq Since $p>n$, for any given $\eg_0>0$ there is $R_0=R_0(\eg_0,C_{1,p})>0$ such that
$$\sup_{R\le R_0}R^{-n}\itQ{Q_R}{|Du|^2}<\eg_0.$$

Since $p>2$, $u$ is bounded by \mref{uniDu} and Sobolev's embedding theorems so that $\llg(u)$ is uniformly bounded, we can apply \cite[Theorem 3.1]{GiaS} to see that $u\in C^{\ag,\ag/2}(Q)$ for all $\ag\in(0,1)$. As $A(\sg u)$ is $C^1$ in $u$, $A(\sg u(x,t))\in C^{\ag,\ag/2}(Q)$ so that \cite[Theorem 3.2]{GiaS} applies to show that $Du\in C^{\ag,\ag/2}_{loc}(Q)$ for all $\ag\in(0,1)$.
In fact, it was shown in \cite{GiaS} that there are $R_*>0$ and a constant $C$ depending on the parameters in A), F) and $C_{1,p}$ such that
\beqno{uniDuzzz}\|Du\|_{C^{\ag,\ag/2}(B_{R}\times(t-R^2,t))} \le C \mbox{ if $R\le R_*$ and } B_{2R}\times(t-4R^2,t)\subset Q.\eeq
The above inequality also holds if the center of $B_R$ is on the boundary $\partial\Og$ so that with $4R^2<t_0$ and $t-4R^2>0$ we see that $\|u\|_{\ccX_3}$ is uniformly bounded. 
 
Thus, we obtain a uniform bound for the fixed points of $T_\sg$ in the Banach space $\ccX$ and the Leray-Schauder theorem can apply to give the existence of a strong solution in $Q$ for any given $T_0>0$. By the uniqueness of strong solutions (as $u$ is bounded and $A,\hat{f}$ are smooth) we see that the strong solutions in two cylinders $Q\subset Q'$ coincide in $Q$. This shows that $u$ is unique and exists globally. We now see that the proof continues to hold without the boundedness assumption of $\llg(x,t,u)$ for $t$ near $0$ in \cite{dleANS}.

Since the systems \mref{e1fam} satisfy A), F) and SG) with the same set of constants so that, in applications, we need only to verify  \mref{Du2bound}, \mref{llgrzzz} and M') for strong solutions to \mref{e1} then the same argument and \reflemm{dleANSprop} shows that the estimate \mref{uniDu} also holds uniformly for the strong solutions of the systems in the family.

The assumption \mref{Du2bound} is usually easy to check by testing the systems with $u$. Meanwhile, \mref{llgrzzz} is also a mild assumption, especially if $\llg(u)$ has polynomial growth as in \mref{e0} and we know that $u$ is BMO (see also  \mref{llgr} and \reflemm{llgs0lem}). Thus, M') is the key assumption needs to be checked in order to establish \mref{uniDu} and then the bound \mref{uniDuzzz} for higher norms of the solution $u$. The numbers $\mu_0, R_{\mu_0}$ are the key parameters determining $C_{1,p}$ and these bounds. One should note that the number $R_{\mu_0}$ may also depend on the initial condition $u(\cdot,0)$ so that the bound $C_{1,p}$ may implicitly depend on $\|u(\cdot,0)\|_{W^{1,p}(\Og)}$.

We are now ready to give the proof of \reftheo{extcoro1} for the planar case $\Og\subset\RR^2$ by checking the conditions of \reftheo{extthm} under the assumption \mref{Du2boundz} that the norm $\|u(\cdot,t)\|_{W^{1,2}(\Og)}$ of any strong solution $u$ of \mref{e1} does not blow up in finite time.

First of all, we observe that the assumption \mref{llgr} in \reftheo{extcoro1} is much weaker than \mref{llgrzzz} in \reflemm{dleANSprop}, which requires higher integrability of $\llg(u)$. This is because we have assumed \mref{Du2boundz} which is a bit stronger than M'). We have the following lemma showing that \mref{llgullgeg},\mref{Du2boundz} and \mref{llgr} imply \mref{llgrzzz}. This lemma will be also used in the proof of  \reftheo{extcoro2}.

\blemm{llgs0lem} Fix a $t\in(0,T_0)$. Assume \mref{llgullgeg} and that there are positives $s_0,M_0, M_1$ and $C_0$ such that \beqno{llgs0}\|\llg^{s_0}(u(\cdot,t))\|_{L^{1}(\Og)}\le C_0M_0^{s_0}, \; \|Du(\cdot,t)\|_{L^2(\Og)}\le M_1. \eeq Then, for any $r>1$ there is a constant $C(C_0,s_0,r,|\Og|,\BLLg_1,M_0,M_1)$ such that \beqno{llgs0r}  \|\llg(u(\cdot,t))\|_{L^{r}(\Og)}\le C(C_0,s_0,r,|\Og|,\BLLg_1,M_0,M_1).\eeq
\elemm

\bproof We choose and fix $s_*>0$ and $p\in(1,2)$ such that $s_*p_*=s_0$, where $p_*=2p/(n-2)$. Then \mref{llgs0} implies \beqno{llgs} \|\llg^{s}(u)\|_{L^{p_*}(\Og)}\le C_0\llg_*,\quad \mbox{ where $s=s_*p_*$ and $\llg_*:=M_0^\frac{s_0}{p_*}$}.\eeq

Fix a $t\in(0,T_0)$ and define $g(\cdot)=\llg^{s+\eg_0}(u(\cdot,t))$. The definition of $\BLLg_1$ in \mref{llgullgeg} gives
$$|Dg|\le C(s)\frac{|\llg_u|}{\llg^{1-\eg_0}(u)}\llg^s(u)|Du| \le C(s)\BLLg_1 \llg^s(u)|Du|.$$ Hence, as $q=(2/p)'$, by H\"older's inequality, $\|Dg\|_{L^p(\Og)} \le C\|\llg^s(u)\|_{L^{p_*}(\Og)}\|Du\|_{L^2(\Og)}$. This implies, using \mref{llgs} and \mref{llgs0}, $\|g\|_{W^{1,p}(\Og)}\le C(C_0,s,\BLLg_1,M_1)\llg_*$. By Sobolev's imbedding theorem, $$\|\llg^{s+\eg_0}(u)\|_{L^{p_*}(\Og)}=\|g\|_{L^{p_*}(\Og)}\le C(C_0,s,\BLLg,M_1)\llg_*.$$  This shows that if \mref{llgs} holds for some $s$ then it also holds for $s$ being $s+\eg_0$ and new constants $C_1$ depending on $M_0,M_1$.  It is clear that we can repeat this argument  to see that $\|\llg^{s+k\eg_0}(u)\|_{L^{p_*}(\Og)}\le C_k\llg_*$ for all integers $k$. A simple use of H\"older's inequality and the definition of $\llg_*$  complete the proof. \eproof

{\bf The proof of \reftheo{extcoro1}} We will make use of \reflemm{dleANSprop} and check its assumptions here. First of all, because we are considering the case $n=2$ the condition SG) is trivially satisfied. Next, thanks to \reflemm{llgs0lem}, it is now clear that the condition \mref{llgrzzz} holds under its weaker version \mref{llgr} and \mref{Du2boundz}.  

Therefore, we only need to show that M') holds for \mref{e1famz}. The argument after \reftheo{extthm}, whose \mref{e1fam} is exactly \mref{e1famz},  then shows that the strong solution $u$ of \mref{e1} exists globally. Since $\BLLg$ is bounded, we need only prove that for some sufficiently small $R$ and any ball $B_R$ and any $t$ in  $(0,T_0)$ the quantity $\|u_\sg(\cdot,t)\|_{BMO(B_{R}\cap\Og)}$ can be arbitrarily small. 
We argue by contradiction.  If this is not the case then there are sequences $\{x_n\}\subset\bar{\Og}$, $\{\sg_n\}\subset [0,1]$, $\{t_n\}\subset(0,T_0)$, $\{r_n\}$, $r_n\to0$, such that  $$\|u_{\sg_n}(\cdot,t_n)\|_{BMO(B_{r_n}(x_n)\cap\Og)}>\eg_0 \mbox{ for some $\eg_0>0$}.$$ 

Let $U_n(\cdot)=u_{\sg_n}(\cdot,t_n)$. By \mref{Du2boundz} we see that the sequence $\{U_n\}$ is bounded in $W^{1,2}(\Og)$. We can then assume that $U_n$ converges weakly to some $U$ in $W^{1,2}(\Og)$ and strongly in $L^2(\Og)$. We then have $\|U_n\|_{BMO(B_R\cap\Og)}\to \|U\|_{BMO(B_R\cap\Og)}$ for any given ball $B_R$.  Since $n=2$, by Poincar\'e's inequality $U$ is VMO and $\|U\|_{BMO(B_R\cap\Og)}<\eg_0/2$ if $R$ is sufficiently small. Furthermore, we can assume also that $x_n$ converges to some $x\in\bar{\Og}$. Thus, for large $n$, we have $r_n<R/2$ and $x_n\in B_{R/2}(x)$. Hence, $B_{r_n}(x_n)\subset B_{R}(x)$ and if $n$ is sufficiently large then $$\|U_n\|_{BMO(B_{r_n}(x_n)\cap\Og)} \le \|U_n\|_{BMO(B_R(x)\cap\Og)}\le \|U\|_{BMO(B_R(x)\cap\Og)}+\eg_0/2<\eg_0.$$ We obtain a contradiction and complete the proof. \eproof

We now give the proof of \reftheo{extcoro2}.

{\bf Proof of \reftheo{extcoro2}} From the theory of global attractors (e.g. see \cite{tem}) it is now well known that we need only establish the following claims:

\bdes \item[Claim 1:] \mref{e1} defines a global semiroup $\{\ccS(t)\}_{t\ge0}$ in the Banach space $X=W^{1,p_0}(\Og)$. Namely, the maps $$\ccS(0)U_0=U_0,\quad \ccS(t)U_0(x)=u(x,t)$$  with $u$ being the solution of \mref{e1} for the given initial data $U_0\in X$, is defined for all $t>0$. \item[Claim 2:] The map $\ccS(t)$ is compact and possesses an absorbing ball in $X$. \edes

As we assume \mref{Du2boundzz} holds for any initial data $U_0\in X$, Claim 1 is already established in \reftheo{extcoro1} proving the global existence of \mref{e1}.

We need only consider Claim 2. 
The discussion after \reftheo{extthm} on the regularity of the solutions shows that $\ccS(t)$ is a compact map on $X$. 
Moreover, \mref{uniDuzzz} shows that the norm $\|u(\cdot,t)\|_{C^{1}(\Og)}$, and thus $\|u(\cdot,t)\|_X$, is bounded in terms of the constant $C_{1,p}$ in \mref{uniDu} of \reflemm{dleANSprop}. We should note that the cutoff function $\eta(t)$ can now be used for some fixed $t_0$ so that the bound $C_{1,p}$ is independent of the initial data $u(\cdot,0)$ but $t_0^{-1}$. Since $u$ exists globally, we can choose $t_0=1$.  Therefore, for any bounded set $K\subset X$ if we can show that the assumption M') is verified uniformly for all $u\in K$ and $t$ sufficiently large then we can choose a universal $C_{1,p}$ and then show that there is an absorbing ball  in $X$.

Thus, we just need to show that the conditions of \reflemm{dleANSprop} hold uniformly for large $t$ to give a uniform bound $C_{1,p}$.  First of all,  for large $t$ in \reflemm{llgs0lem} we let $M_0=M^{1/s_0}$ and $M_1=M$, where $M$ is the universal constant in \mref{Du2boundzz} and \mref{llgrzz}. By \reflemm{llgs0lem},
 the condition \mref{llgrzzz} of \reflemm{dleANSprop} holds uniformly when $t$ is large.

To prove that M') holds uniformly, we use a contradiction argument similar to that in the proof of \reftheo{extcoro1}. We need only show that for some sufficiently small $R$ and any ball $B_R$ and any bounded set $K$  the quantity $\|u(\cdot,t)\|_{BMO(B_{R}\cap\Og)}$, where $u$ is a strong solution with initial data $U_0$ in $K$, can be arbitrarily small when $t$ is sufficiently large.  If this is not the case then there are sequences $\{K_n\}$ of bounded sets  in $X$, $\{u_n\}$ of strong solutions to \mref{e1} with initial data in $K_n$, $\{x_n\}\subset\bar{\Og}$, $\{t_n\}\subset(0,\infty)$ with $t_n>T_{K_n}$, $\{r_n\}$, $r_n\to0$, such that  $$\|u_n(\cdot,t_n)\|_{BMO(B_{r_n}(x_n)\cap\Og)}>\eg_0 \mbox{ for some $\eg_0>0$}.$$ 

Let $U_n(\cdot)=u_n(\cdot,t_n)$. By the same argument in the end of the proof of \reftheo{extcoro1}  we obtain a contradiction and complete the proof. \eproof

\section{The Generalized SKT systems on Planar Domains}\eqnoset\label{2dsktsec}

In this section, we consider the generalized version of the SKT system \mref{e0}
\beqno{e00}u_t=\Delta(P(u))+\hat{f}(u,Du)\quad  (x,t)\in Q=\Og\times(0,T),\eeq
where $u:\Og\to\RR^m$, $P:\RR^m\to\RR^m$ and $f:\RR^m\times\RR^{nm}\to\RR^m$ are vector valued functions.  It is clear that \mref{e00} generalizes \mref{e0}, where $P(u)$ is quadratic in $u$, but is still a special case of \mref{e1} for  $A(u)=P_u(u)$
$$u_t=\Div(A(u)Du)+\hat{f}(u,Du)\quad (x,t)\in Q=\Og\times(0,T_0).$$

Inspired by the usual SKT system \mref{e0}, we will allow the following polynomial growth condition of $A(u), \hat{f}(u)$.

\bdes

\item[SKT)] Assume that $P(0)=0$ and $A(u):=P_u(u)$ satisfies A). Moreover, there are positive constants $\llg_0,\llg_1,k,C$ and $r_0$ such that  for all $u\in\RR^m$    \beqno{Fghyp} \llg(u)\sim \llg_0+\llg_1|u|^k,\quad |\llg_u(u)| \le \left\{\barr{l}\mbox{bounded by $C\llg_1|u|^{k-1}$ if $|u|>r_0$},\\\mbox{$C$ if $|u|\le r_0$.}\earr\right.\eeq  
In addition,  $|P(u)|\le C\llg(u)|u|$.
\edes

Here and in the sequel, we will write $a\sim b$ if there are two generic positive constants $C_1,C_2$ such that $C_1b \le a \le C_2b$.

Firstly, it is easy to check that the assumption \mref{LLgx} in A) on the finiteness of $$ \BLLg=\sup_{u\in\RR^m}\frac{|\llg_u(u)|}{\llg(u)}$$is satisfied under the above polynomial growth assumption \mref{Fghyp} on $\llg,\llg_u$.

Secondly, since $P(0)=0$ and $A(u)=P_u(u)$ with  $|A(u)|\sim \llg(u)$, it is natural to assume $|P(u)|\le C\llg(u)|u|$ for some constant $C$. 

We also note that $\llg(u)$ is the smallest eigenvalue of $(A+A^T)/2$ and $\llg^2(u)$ is the smallest eigenvalue of $A^TA$. Therefore, \mref{A1} implies
\beqno{A2}C\llg^2(u)|\zeta|^2\ge |A(u)\zeta|^2=\myprod{A^T(u)A(u)\zeta,\zeta}\ge  \llg^2(u)|\zeta|^2.\eeq

Similarly, we will assume the following assumption on the reaction terms.

\bdes\item[F')] $\hat{f}$ satisfies F) and there is $C_f>0$ such that for $\llg_S=\llg_0+\llg_1$
\beqno{fllg}|f(u)|\le C_f\llg_S^{-1}|u|\llg(u)\quad \forall u\in\RR^m.\eeq
\edes

We note that if $f(0)=0$ then \mref{fUDU21} in F) implies  \mref{fllg} for some suitable constant $C_f$.

Our first result in this section shows that global existence of \mref{e00} can be established under much weaker assumption on the $L^p$ norm of strong solutions to the following systems.
\beqno{e1famzz}\left\{\barr{l} u_t=\Div(A(\sg u)Du)+\hat{f}(\sg u,\sg Du)\quad (x,t)\in Q=\Og\times(0,T_0), \sg\in[0,1]\\u(x,0)=U_0(x)\quad  x\in\Og\\\mbox{$u=0$ or $\frac{\partial u}{\partial \nu}=0$ on $\partial \Og\times(0,T_0)$}. \earr\right.\eeq

\btheo{thmglobext} Assume SKT), F'). Let $u_\sg$ be a strong solution to \mref{e1famzz} and $(0,T_0)$ be its maximal existence interval. Suppose that there are positive constants  $q>1$ and $M_{U_0,T_0}$ depending on $\|U_0\|_{W^{1,p_0}(\Og)}$ and $T_0$  such that \beqno{2kboundz}\|u_\sg(\cdot,t)\|_{L^{qk}(\Og)}\le M_{U_0,T_0} \quad \mbox{for all $t\in (0,T_0)$}.\eeq 

Then 
\beqno{DuTest0}\iidx{\Og\times\{t\}}{\llg(u_\sg)|Du_\sg|^2}\le C(T_0,C_f, M_{U_0,T_0})\llg_S\llg_0^{-1} \quad \forall t\in(0,T_0).\eeq

Moreover, \mref{e00} or \mref{e1famzz} for $\sg=1$ has a unique strong solution $u$ which exists globally, i.e. $T_0=\infty$.

\etheo

Next, if we can control the bound \mref{2kboundz} uniformly when $t$ is large then we have the following result on the existence of global attractors.

\btheo{thmglobatt} Assume SKT), F'). Suppose that there are constants $q>1$ and $M$ such that the dynamical system defined by \mref{e00} possesses an absorbing ball in $L^{qk}(\Og)$. That is there is a constant $M$ such that for any bounded set $K\subset W^{1,p_0}(\Og)$ there is a $T_K>0$ and any global solution $u$ of \mref{e00} with initial data in $K$ will satisfy  \beqno{2kboundzz}\|u(\cdot,t)\|_{L^{qk}(\Og)}\le M \quad \mbox{for all $t\ge T_K$}.\eeq 

Then 
\beqno{DuTest00}\iidx{\Og\times\{t\}}{\llg(u)|Du|^2}\le C(M)\llg_S\llg_0^{-1} \quad \mbox{for all $t\ge T_K+1$}.\eeq

Moreover, the system \mref{e00} possesses a global attractor in $X=W^{1,p_0}(\Og)$.

\etheo

The proof of the above theorems will be based on the checking of the assumptions of on the boundedness of $\|u\|_{W^{1,2}(\Og)}$ and $\|\llg^{s_0}(u)\|_{L^{1}(\Og)}$ in \reftheo{extcoro1} and \reftheo{extcoro2}. By Sobolev's embedding theorems for $n=2$ and the polynomial growth assumption on $\llg(u)$, we need only establish the corresponding boundedness of $\|u\|_{W^{1,2}(\Og)}$. This is clear from \mref{DuTest0}, \mref{DuTest00} and the fact that $\llg(u)$ is bounded from below.

In the sequel,  when there is no ambiguity $C, C_i$ will denote universal constants that can change from line to line in our argument. If necessary, $C(\cdots)$ is used to denote quantities which are bounded in terms of theirs parameters.  Furthermore, as we will always consider a strong solution to \mref{e00} that exists in its maximal time interval $(0,T_0)$, the derivatives $u_t, D^2u$ and $Du_t$ make sense in the proof below.

We first have the following lemma establishing a differential inequality which will be used in several places later on.

\blemm{ldulemm0z} Assume A), F) and that $A(u)=P_u(u)$ for some $P:\RR^m\to\RR^m$. Let $u$ be a strong solution to \mref{e00} on some interval $(0,T_0)$. For any nonnegative $C^1$ function $\eta(t)$ on $[0,\infty)$ we have \beqno{Aeqn}\barr{ll}\lefteqn{\iidx{\Og}{\llg(u)|u_t|^2\eta^2} +\frac{d }{d t}\iidx{\Og}{|A(u)Du|^2\eta^2}\le}\hspace{2cm}&\\& C\iidx{\Og}{[|A(u)Du|^2(\eta\eta_t+1)+\llg(u)|f(u)|^2\eta^2]}.\earr\eeq

\elemm

\bproof 
We test the system of $u$ by $A(u) u_t\eta^2(t)$ (i.e. multiplying the $i^{th}$ equation of \mref{e00} by $\sum_j a_{ij}(u)(u_j)_t\eta^2$, $A(u)=(a_{ij}(u))$, integrating over $\Og$ and summing the results) and integrate by parts in $x$ to get
$$\iidx{\Og}{(\myprod{A(u)u_t,u_t} +\myprod{A(u)Du,D(A(u)u_t)})\eta^2}=\iidx{\Og}{\myprod{\hat{f}(u,Du),A(u)u_t}\eta^2}.$$

As $D(A(u)u_t)=D(P(u)_t)=(DP(u))_t=(A(u)Du)_t$, we see that $$\frac12\frac{\partial }{\partial t}(|ADu|^2\eta^2)=\myprod{A(u)Du,D(A(u)u_t)}\eta^2+|A(u)Du|^2\eta\eta_t,$$ and obtain
\beqno{ztemp}\barr{cc}\lefteqn{\iidx{\Og}{[\myprod{A(u)u_t,u_t}\eta^2 +\frac12\frac{\partial }{\partial t}(|ADu|^2\eta^2)]}=}\hspace{3cm}&\\&\iidx{\Og}{[|A(u)Du|^2\eta\eta_t+\myprod{\hat{f}(u,Du),A(u)u_t}\eta^2]}.\earr\eeq

We now use the ellipticity of $A(u)$ in the fisrt integrand on the left hand side to have $\myprod{A(u)u_t,u_t}\ge\llg(u)|u_t|^2$. Also, as $|\hat{f}(u,Du)| \le C\llg^\frac12(u)|Du| + f(u)$ and $|A(u)|\le C\llg(u)$, we use Young's inequality to find a constant $C(\eg)$ such that for any $\eg>0$ we can estimate the second integrand on the right hand side as follows
$$|\myprod{\hat{f}(u,Du),A(u)u_t}|\le  \eg\llg(u)|u_t|^2 + C(\eg)[\llg^2(u)|Du|^2+\llg(u)|f(u)|^2].$$
Hence, using this in \mref{ztemp} with  sufficiently small $\eg$ and noting that $|A(u)Du|^2\sim\llg^2(u)|Du|^2$, see \mref{A2}, we get \mref{Aeqn}. \eproof 

The integrand $\llg(u)|f(u)|^2\eta^2$ on the right hand side of \mref{Aeqn} will play an important role in our analysis so that the following lemma will show that it can be controlled under some boundedness assumption on the $L^p$ norm of $u$.

To proceed we collect some well known inequalties here for later use. For any $p\ge1$, $\ag\in(0,1)$ and $\eg>0$ we have the following inequality for all $w\in W^{1,2}(\Og)$, recalling that $n=2$ 
\beqno{poinc} \left(\iidx{\Og}{w^p}\right)^\frac2p \le \eg\iidx{\Og}{|Dw|^2} +C(\eg)\left(\iidx{\Og}{w^\ag}\right)^\frac2\ag.\eeq

Concerning the last integral, for $w=\llg_0|u|^p+\llg_1|u|^q$ and suficiently small  $\ag$ we note the following simple fact which results from  H\"older's inequality
\beqno{llgP}\left(\iidx{\Og}{(\llg_0|u|^p+\llg_1|u|^q)^\ag}\right)^\frac2\ag\le C(\ag)(\llg_0^2\|u\|_{L^1(\Og)}^{2p}+\llg_1^2\|u\|_{L^1(\Og)}^{2q}).\eeq

Combining the above two inequalities, for $w=|u|^{\frac q2+1}$, $\llg_0=\llg_1=1$ and $p=2$, we have
\beqno{poinc1} \iidx{\Og}{|u|^{q+2}} \le \eg\iidx{\Og}{|u|^q|Du|^2} +C(\eg,q)\|u\|_{L^1(\Og)}^{q+2}.\eeq

Similarly, letting $w=\llg(u)|u|$  and noting that $w\sim |u|^{k+1}$  and  $|Dw|\sim \llg(u)|Du|$, for any $p\ge1$ we can find a constant $C_p(\eg,\|u\|_{L^1(\Og)})$ such that (with $\llg_S=\llg_0+\llg_1$)
\beqno{poinc1z}\left(\iidx{\Og}{(\llg(u)|u|)^{p}}\right)^\frac2{p} \le \eg\iidx{\Og}{\llg^2(u)|Du|^2} + C_p(\eg,\|u\|_{L^1(\Og)})\llg_S^2.\eeq

\blemm{ldulemm00z}  For some $q>1$ and $t\in(0,T_0)$ we suppose that the number $M(t)=\|u(\cdot,t)\|_{L^{qk}(\Og)}$ is finite. Then, for any $\eg_0>0$ we can find a positive constant $C(\eg_0,M(t))$ such that
\beqno{llgfest}\iidx{\Og\times\{t\}}{\llg(u)|f(u)|^2}\le  C_f^2\llg_S^{-1}\left[\eg_0\iidx{\Og\times\{t\}}{|A(u)Du|^2}+C(\eg_0,M(t))\llg_S^2\right].\eeq

\elemm

\bproof We write \beqno{llgf2z}\llg(u)|f(u)|^2= \llg^2(u)|u|^2\LLg(u), \mbox{ where } \LLg(u) = \frac{|f(u)|^2}{|u|^2\llg(u)}.\eeq For any $q>1$, and $q'=q/(q-1)$ we can use  H\"older's inequality to have
$$\iidx{\Og}{\llg(u)|f(u)|^2}= \iidx{\Og}{\llg^2(u)|u|^2\LLg(u)}\le  C\left(\iidx{\Og}{(\llg(u)|u|)^{2q'}}\right)^\frac1{q'}
\left(\iidx{\Og}{\LLg^{q}(u)}\right)^\frac1{q}.$$

For the first factor on the right we use \mref{poinc1z} with $w=\llg(u)|u|$ with $p=2q'$ and the fact that $\|u(\cdot,t)\|_{L^{1}(\Og)}$ is bounded by $M(t)=\|u(\cdot,t)\|_{L^{qk}(\Og)}$ to find a constant $C_0(\eg,M(t))$ such that
$$\left(\iidx{\Og}{(\llg(u)|u|)^{2q'}}\right)^\frac1{q'} \le \eg\iidx{\Og}{\llg^2(u)|Du|^2} + C_0(\eg,M(t))\llg_S^2.$$ 

On the other hand, from the assumption \mref{fllg} in F'), for some constant $C$ we have $\LLg(u)\le C_f^2\llg_S^{-2}\llg(u)$
so that, using the growth condition on $\llg(u)$ $$\left(\iidx{\Og}{\LLg^{q}(u)}\right)^\frac1{q}\le C_f^2\llg_S^{-2}\|\llg(u)\|_{L^{q}(\Og)}\le C_f^2\llg_S^{-1}(1+\|u\|_{L^{kq}(\Og)})\le C_f^2C_1(M(t))\llg_S^{-1}.$$ Therefore, for any given $\eg_0>0$ we choose $\eg=\eg_0C_1^{-1}(M(t))$ and combine the above estimates to obtain a constant $C(\eg_0,M(t))$ such that
$$\iidx{\Og}{\llg(u)|f(u)|^2}\le  C_f^2\llg_S^{-1}\left[\eg_0\iidx{\Og}{\llg^2(u)|Du|^2}+C(\eg_0,M(t))\llg_S^2\right].$$ By \mref{A2}, $|A(u)Du|^2\sim \llg^2(u)|Du|^2$, the above proves the lemma.
\eproof

We now apply \reftheo{extcoro1} to establish the global existence result.

{\bf Proof of \reftheo{thmglobext}:} 
From \mref{Aeqn} of \reflemm{ldulemm0z} with $\eta\equiv1$, we see that
\beqno{p2}\frac{d }{d t}\iidx{\Og}{|ADu|^2}\le C\iidx{\Og}{(|A(u)Du|^2+\llg(u)|f(u)|^2)}.
\eeq

We let $\eg_0=1$ in \mref{llgfest} and replace $M(t)$ by $M_{U_0,T_0}$, see  \mref{2kboundz}, to obtain a constant  $C(M_{U_0,T_0})$ such that $$\iidx{\Og}{\llg(u)|f(u)|^2}\le  C_f^2\llg_S^{-1}\left[\iidx{\Og}{|A(u)Du|^2}+C(M_{U_0,T_0})\llg_S^2\right].$$

We then have for $\bg:=C_f^2\llg_S^{-1}$ and some constant $C$ the following inequality. \beqno{p3azzzz}\frac{d }{d t}\iidx{\Og}{|ADu|^2}\le C(\bg+1)\iidx{\Og}{|A(u)Du|^2}+C_f^2C(M_{U_0,T_0})\llg_S.
\eeq 

We now set  $$y(t)=\iidx{\Og}{|A(u)Du(x,t)|^2}$$ to see from \mref{p3azzzz} that $y'(t) \le C(\bg+1) y(t)+ C(M_{U_0,T_0},C_f)\llg_S$ for all $t\in(0,T_0)$. A simple use of Gronwall's inequality  shows that there exists a constant $C(T_0, C_f,M_{U_0,T_0})$ such that
\beqno{DuTest}\iidx{\Og\times\{t\}}{|A(u)Du|^2}\le C(T_0,C_f, M_{U_0,T_0})\llg_S \quad \forall t\in(0,T_0).\eeq
 
 Because $|A(u)Du|^2\ge C\llg^2(u)|Du|^2$ for some $C>0$ and $\llg(u)$ is bounded from below by $\llg_0$, the above yields the bound \mref{DuTest0} and a bound for $\|Du(\cdot,t)\|_{L^2(\Og)}$ for all $t>0$. We see that \reftheo{extcoro1} applies here. In fact, it is easy to check that the data of the systems in \mref{e1famz} satisfies SKT) and F') with the same set of constants so that the above argument (for $\sg=1$) yields a uniform estimate for $\|u_\sg(\cdot,t)\|_{W^{1,2}(\Og)}$ in \mref{Du2boundz} of \reftheo{extcoro1}. The bound for  $\|\llg^{s_0}(u_\sg(\cdot,t))\|_{L^1(\Og)}$ in \mref{llgr} follows from the polynomial growth of $\llg(u)$ and Sobolev's inequality. The proof is then complete. \eproof

We now turn to the proof of \reftheo{thmglobatt} and apply \reftheo{extcoro2}. We need only show that there is an absorbing ball in $W^{1,2}(\Og)$. To proceed we need some lemmas establishing uniform estimates for $\|Du\|_{L^2(\Og)}$ under much weaker assumptions on the $L^p$ norms of $u$.

\blemm{ldulemm1z}  Let $u$ be a strong solution to \mref{e00} on some interval $(0,T_0)$. For any $\tau_0\in(0,1)$ assume that there are positive constants  $q>1$,  $T\ge0$ and $T'=\min\{T+\tau_0,T_0\}$ such that  the number \beqno{2kbound}M_{T,T',q}=\sup_{t\in[T,T']}\|u(\cdot,t)\|_{L^{qk}(\Og)}\eeq is finite. Then there exists a positive constant $C_0(M_{T,T',q},C_f)$ such that
\beqno{DuTestz}\dspl{\int}_{T+\tau_0}^{T'}\iidx{\Og}{|A(u)Du|^2}ds \le C_0(M_{T,T',q},C_f)(T'-T)(\frac{1}{\tau_0}+\llg_S).\eeq \elemm

\bproof  For $T\ge0$, $\tau_0>0$ and $T'=\min\{T+\tau_0,T_0\}$ and any $R\in(0,\tau_0)$ let  us denote  $$I_R=[T+\tau_0-R,T'].$$ Also, for any $0\le\rg<R\le \tau_0$ we let $\eta$ be a cut-off function for $I_R,I_\rg$. That is, $\eta\equiv 1$ in $I_\rg$ and $\eta(\tau)\equiv 0$ for $\tau\le T+\tau_0-R$ and $|\eta_t|\le 1/(R-\rg)$.

We multiply the system \mref{e00} by $P(u)\eta^2(t)$  and integrate by parts in $x$ to get
$$\iidx{\Og}{\myprod{A(u)Du,A(u)Du}\eta^2}=\iidx{\Og}{(-\myprod{P(u),u_t} +\myprod{\hat{f}(u,Du),P(u)})\eta^2}.$$

Since $|\hat{f}(u,Du)| \le C\llg^\frac12(u)|Du| + f(u)$ and $|A(u)Du|^2\sim\llg^2(u)|Du|^2$ (see \mref{A2}), we apply Young's inequality to the first term in the integrand on the right hand side and find a constant $C_1$ such that for any $\eg_*>0$
$$\iidx{\Og}{|A(u)Du|^2\eta^2} \le C_1\iidx{\Og}{(\eg_*\llg(u)|u_t|^2+\eg_*^{-1}\llg^{-1}(u)|P(u)|^2+|f(u)||P(u)|)\eta^2}.$$
Integrate the above over  $I_R$ to get
\beqno{yz1}\barr{ll}\lefteqn{\dspl{\int_{I_R}}\iidx{\Og}{|A(u)Du|^2\eta^2}ds \le}\hspace{2cm}&\\& C_1\dspl{\int_{I_R}}\iidx{\Og}{(\eg_*\llg(u)|u_t|^2+\eg_*^{-1}\llg^{-1}(u)|P(u)|^2+|f(u)||P(u)|)\eta^2}ds.\earr\eeq

We integrate \mref{Aeqn} of \reflemm{ldulemm0z} over $I_R$ and note that $|ADu|^2\eta^2=0$ at $T+\tau_0-R$. From the choice of $\eta$ we can find a constant $C_2$ to obtain
$$\barr{cc}\lefteqn{\int_{I_R}\iidx{\Og}{\llg(u)|u_t|^2\eta^2}ds +\iidx{\Og\times\{T'\}}{|A(u)Du|^2}\le}\hspace{3cm}&\\& C_2\dspl{\int}_{I_R}\iidx{\Og}{[(\frac{1}{R-\rg }+1)|A(u)Du|^2+\llg(u)|f(u)|^2\eta^2]}ds.\earr$$
Hence, \beqno{yz0}\int_{I_R}\iidx{\Og}{\llg(u)|u_t|^2\eta^2}ds\le C_2\int_{I_R}\iidx{\Og}{[(\frac{1}{R-\rg }+1)|A(u)Du|^2+\llg(u)|f(u)|^2\eta^2]}ds.\eeq

We now take $\eg_*=\frac12(C_1C_2)^{-1}(R-\rg)$ in \mref{yz1} and note that  $\eg_*\le C\tau_0$ for some fixed constant $C$. Hence, multiplying \mref{yz0} by $C_1\eg_*$ and using the result in \mref{yz1} we easily get (using $I_\rg \subset I_R$)
\beqno{yz2}\barr{ll}\lefteqn{\dspl{\int}_{I_\rg}\iidx{\Og}{|A(u)Du|^2}ds \le \frac12(1+R-\rg )\dspl{\int}_{I_R}\iidx{\Og}{|A(u)Du|^2}ds+}\hspace{.5cm}&\\&
C_3\dspl{\int}_{I_R}\iidx{\Og}{\left(\frac{1}{|R-\rg|}\llg^{-1}(u)|P(u)|^2+|f(u)||P(u)|+\llg(u)|f(u)|^2\right)}ds.\earr\eeq

We then estimate the last integral on the right hand side of \mref{yz2}. In the sequel, we will abbreviate $$M:=M_{T,T',q} \mbox{ and } |I|=T'-T.$$

Firstly,  from SKT) we have $\llg^{-1}(u)|P(u)|^2\le C\llg(u)|u|^2$. The integral of $\llg(u)|u|^2$ can be estimated by using \mref{poinc} for $w=\llg^\frac12(u)u$ and the fact that $\|u(\cdot,t)\|_{L^1(\Og)}$ is bounded in terms of $M$ for $t\in (T,T')$. We then have
$$\barr{lll}\dspl{\int}_{I_R}\iidx{\Og}{\llg^{-1}(u)|P(u)|^2}ds &\le& C\dspl{\int}_{I_R}\iidx{\Og}{(\llg(u)|Du|^2+C_4(M)\llg_S)}ds\\&\le& C_5(M)|I|(1+\llg_S).\earr$$ Here, we have used the fact that there is a constant $C(M,\tau_0)$ such that \beqno{llgDu2z}\int_{I_R}\iidx{\Og}{\llg(u)|Du|^2} \le C(M,\tau_0)|I|.\eeq This can be proved easily by testing the system with $u$ (see also \refrem{fduremz} for details).

Next, using \mref{poinc} for $w=P(u)$ and any $\eg_0>0$ we can find $C(\eg_0,M)$ such that
\beqno{Pest}\dspl{\int}_{I_R}\iidx{\Og}{|P(u)|^2}ds \le \dspl{\int}_{I_R}\iidx{\Og}{(\eg_0|A(u)Du|^2+C(\eg_0,M)\llg_S^2)}ds.\eeq By F'), $|f(u)|\le C_f\llg_S^{-1}|u|\llg(u)$. Applying Young's inequality, we have
$$|f(u)||P(u)|\le C_f^{-1}\llg_S|f(u)|^2 + C_f\llg_S^{-1}|P(u)|^2\le C_f\llg_S^{-1}(|u|^2|\llg(u)|^2 + |P(u)|^2).$$
Using \mref{poinc1z} ($p=2$) and \mref{Pest} to estimate the integral of the right hand side of  the above, we obtain 
$$\dspl{\int}_{I_R}\iidx{\Og}{|f(u)||P(u)|}ds \le C_f\llg_S^{-1}\dspl{\int}_{I_R}\iidx{\Og}{(\eg_0|A(u)Du|^2+C(\eg_0,M)\llg_S^2)}ds.$$

Concerning the last integrand on the right hand side of \mref{yz2}, by \mref{llgfest}, we have $|f(u)|\le C_f\llg_S^{-1}|u|\llg(u)$
$$\iidx{\Og}{\llg(u)|f(u)|^2}\le  C_f^2\llg_S^{-1}\left[\eg_0\iidx{\Og}{|A(u)Du|^2}+C(\eg_0,M)\llg_S^2\right].$$

We now set for $\tau\in[0,\tau_0]$ $$F(\tau):=\dspl{\int}_{I_\tau}\iidx{\Og}{|A(u)Du|^2}ds$$ and use the above estimates for the integrals on the right hand side of \mref{yz2} to obtain \beqno{Fiter}F(\rg) \le \mu F(R)+\frac{\llg_SC_6(M)}{R-\rg}|I| + C(\eg_0,M,C_f)|I|(1+\llg_S),\eeq where for some constant $C_7$ \beqno{mudef}\mu=\frac12(1+R-\rg )+C_7(C_f+C_f^2)\eg_0\llg_S^{-1}.\eeq

Since $\llg_S$ is bounded from below and $R-\rg\le\tau_0<1$, we can choose $\eg_0$ sufficiently small and depending only on $C_f$ such that $\mu<1$. By an elementary iteration lemma \cite[Lemma 6.1, p.192]{Gius}, we have
$$F(\rg) \le C\left[\frac{C_6(M)\llg_S}{R-\rg}|I| + C_8(M,C_f)|I|\llg_S\right] \mbox{ for all $0\le \rg<R\le\tau_0$}.$$

We now take $R=\tau_0$, $\rg =0$ in the above inequality to obtain
$$\dspl{\int}_{T+\tau_0}^{T'}\iidx{\Og}{|A(u)Du|^2}ds \le C(M,C_f)|I|\llg_S(\frac{1}{\tau_0}+1).$$
The above estimate completes the proof as $M:=M_{T,T',q}$ and $|I|=T'-T$. \eproof

We are now ready to give

{\bf The proof of \reftheo{thmglobatt}:} We will apply \reftheo{extcoro2} here. To this end, let $K$ be any bounded subset of $X$ and $u$ be the solution of \mref{e00} with initial data $U_0$ in $K$.  Under the assumption \mref{2kboundzz} \beqno{mtu}\|u(\cdot,t)\|_{L^{qk}(\Og)}\le M\quad \forall t\ge T_K,\eeq we will show that there exists a constant $C(M)$  such that
\beqno{DuTest}\iidx{\Og\times\{t\}}{\llg^2(u)|Du|^2}\le C(M)\llg_S \quad \forall t\ge T_K+1.\eeq

Because $\llg(u)$ is bounded from below by $\llg_0$, \mref{DuTest} yields a uniform bound for $\|Du(\cdot,t)\|_{L^2(\Og)}$, $t>T_K+1$. \reftheo{extcoro2} then applies and completes the proof.

First of all, by the global existence result we can assume $T_K>1$ and take $\tau_0=1/2$ in \reflemm{ldulemm1z}.
From \mref{llgfest} with $\eg_0=1$ and suitable choice of $C(M)$, we have $$\iidx{\Og}{\llg(u)|f(u)|^2}\le  C_f^2\llg_S^{-1}\left[\iidx{\Og}{|A(u)Du|^2}+C(M)\llg_S^2\right]\quad \forall t\ge T_K.$$

As in the proof of \reftheo{thmglobext}, see \mref{p3azzzz}, we obtain for $\bg:=C_f^2\llg_S^{-1}$ and $\ag:=C_f^2C(M)$ \beqno{p3az}\frac{d }{d t}\iidx{\Og}{|A(u)Du|^2}\le C(\bg+1)\iidx{\Og}{|A(u)Du|^2}+\ag\llg_S \quad \forall t\ge T_K.
\eeq

We now set  $$y(t)=\llg_S^{-1}\iidx{\Og}{|A(u)Du(x,t)|^2}$$ and derive from \mref{p3az} that $$y'(t) \le C(\bg+1) y(t) +\ag\quad \forall t\in(T_K,\infty).$$

For any $\tau_0>0$ the uniform Gronwall inequality (see \cite[Lemma 1.1, p.91]{tem}) then gives
\beqno{unigronwall}y(t+\tau_0) \le\left[\frac{a_3}{\tau_0}+a_2\right]\exp(a_1) \quad \forall t>T_K,\eeq where $$a_1:=\int_t^{t+\tau_0}C(\bg+1) ds,\, a_2:=\int_t^{t+\tau_0}\ag ds,\,a_3:=\int_t^{t+\tau_0}y(s)ds.$$

From the definitions of $\ag,\bg$
$$a_1\le C(C(M)C_f^2\llg_S^{-1}+1)\tau_0 \mbox{ and }a_2\le C(M)C_f^2\tau_0.$$

For a fixed $\tau_0\in(0,1)$ and any $t\ge T_K+\tau_0$ we take $T$ and $T'$ in \mref{DuTestz} such that $T+\tau_0=t$ and $T'=t+\tau_0$. Then $T\ge T_K$ so that  \reflemm{ldulemm1z} applies here to give 
$$\frac{a_3}{\tau_0}=(\tau_0\llg_S)^{-1}\int_t^{t+\tau_0}\iidx{\Og}{|A(u)Du(x,t)|^2}ds \le\llg_S^{-1} C_0(M,C_f)(\frac{1}{\tau_0}+\llg_S).$$

 Putting these estimates in \mref{unigronwall} and the fact that $\llg_S$ is bounded from below, we see that $y(t+\tau_0)$ is uniformly bounded by a constant $C(M)$ for $t>T_K+\tau_0$. Using the definition of $y(t)$ and letting $\tau_0=1/2$ we obtain \mref{DuTest} and complete the proof. \eproof

Next, we will show that the assumption on the boundedness of $L^p$ norm of the solutions in \mref{2kboundz} and \mref{2kboundzz} can be weaken further if a mild hyposthesis on the structure on $A(u)$ is imposed. In fact, the $L^p$ norm can be replaced by $L^1$ one if we assume further the following (the constant $C_*$ is the ratio between the eigenvalues of $A(u)$ described in A)).
\bdes\item[SG')] If $k>2$ then there is a number $\dg_k\in(0,1)$ such that $(k-2)/k\le \dg_k C_*^{-1}$. \edes

\btheo{L1thm} Assume SKT), F') and SG'). The conclusions of \reftheo{thmglobext} and \reftheo{thmglobatt} hold if  \mref{2kboundz} and \mref{2kboundzz} respectively hold for the $L^1(\Og)$ norm of $u$. \etheo

The proof of this theorem clearly follows from the following lemma (and \refrem{tau0is0}) which shows that an appropriate bound for the $L^1$ norm of $u$ implies those of $L^{qk}$ norm of $u$ for some $q>1$. \reftheo{thmglobext} and \reftheo{thmglobatt} then apply.

\blemm{ldulemm} Assume SKG), F') and SG'). Suppose that there are $T_*\in(0,T_0)$ and a continuous function $C_0$ on $(T_*,\infty)$ such that \beqno{u1uni} \|u(\cdot,t\|_{L^1(\Og)} \le C_0(t)\quad\mbox{for all $t\in(T_*,T_0)$}.\eeq

Then, for any positive $\tau_0<T_0-T_*$ there is  a number $q>1$ such that
\beqno{p00z} \iidx{\Og}{|u(x,t)|^{qk}} \le C(\sup_{(T_*,T_0)}C_0(t),\tau_0) \quad \forall t\in(T_*+\tau_0,T_0).\eeq
\elemm

\bproof  First of all, we recall the following fact from \cite{letrans} if $l>0$ and $$\frac{l}{l+2}\le \dg_l C_*^{-1} \mbox{ for some } \dg_l\in(0,1)$$ then there exists a positive $\llg_l$, which depends on $k$, such that \beqno{sgineq}\myprod{A(u)Du,D(|u|^l u)}\ge \llg_l\llg(u) |u|^l|Du|^2.\eeq
By SG') it is clear that we can find $l>\max\{0,k-2\}$ such that \mref{sgineq} holds. 

Let $T\ge T_*$ such that there is $T'\in(T+\tau_0,T_0)$. We test the system of $u$ with $|u|^lu\eta^p(t)$, where $\eta$ is a cutoff function for $[T,T'], [T+\tau_0,T']$  and $p>1$, which will be determined shortly. Using \mref{sgineq}  and the fact that $|\eta_t|\le 1/\tau_0$, we easily obtain for $Q=\Og\times[T,T']$
\beqno{p0}\barr{ll}\lefteqn{\sup_{t\in [T+\tau_0,T']}\iidx{\Og}{|u|^{l+2}} + \llg_l\itQ{Q}{\llg(u) |u|^l|Du|^2\eta^p} \le}\hspace{5cm}&\\& C\itQ{Q}{[\myprod{\hat{f}(u,Du),|u|^lu}\eta^p+\frac{1}{\tau_0}|u|^{l+2}\eta^{p-1}]}.\earr\eeq

For the last term in the integrant on the right, we can choose $p$ such that  $p-1 > p\frac{l+2}{k+l+2}$ and use Young's inequality to find some positive constant $C(k,\tau_0)$ such that
$$\frac{1}{\tau_0}|u|^{l+2}\eta^{p-1} \le |u|^{k+l+2}\eta^p + C(k,\tau_0)\le |u|^{k+l+2}\eta^p + C(k,\tau_0).$$

By SKT) and F'), we can use Young's inequality to find a constant $C$ such that $$\myprod{\hat{f}(u,Du),|u|^lu} \le \eg\llg(u)|u^l|Du|^2 + C(\eg)(|u|^{k+l+2}+1).$$  Using \mref{poinc1} for $q=k+l$ and the assumption \mref{u1uni}, we have for $t\in(T_*,T')$ $$\barr{lll}\iidx{\Og\times\{t\}}{|u|^{k+l+2}\eta^p} &\le& \eg\iidx{\Og\times\{t\}}{|u|^{k+l}|Du|^2\eta^p}+C(\eg,C_0(t))\\&\le& \eg\iidx{\Og\times\{t\}}{|u|^{l}\llg(u)|Du|^2\eta^p}+C(\eg,C_0(t)).\earr$$

Thus, there is a constant $C(\eg,\sup_{(T,T')}C_0(t),k,\tau_0)$ such that the right hand side of \mref{p0} can be estimated by $$\eg\itQ{Q}{|u|^{l}\llg(u)|Du|^2\eta^p}+C(\eg,\sup_{(T,T')}C_0(t),k,\tau_0)|\Og|(T'-T).$$ Choosing $\eg=\llg_l/2$, we can deduce from \mref{p0} and the above inequality the following estimate. 
\beqno{llgDu2}\sup_{t\in [T+\tau_0,T']}\iidx{\Og}{|u|^{l+2}} + \itQ{Q}{|u|^{l}\llg(u)|Du|^2\eta^p} \le C(\sup_{(T,T')}C_0(t),k,\tau_0)|\Og|(T'-T).\eeq Therefore, $$\iidx{\Og}{|u(x,t)|^{l+2}} \le C(\sup_{(T_*,T')}C_0(t),k,\tau_0,T_0)|\Og| \quad \forall t\ge T+\tau_0.$$

Since $l+2>k$, there is $q>1$ such that $l+2=qk$ and the above yields \mref{p00z}. The proof is then complete. \eproof

\brem{tau0is0} We can allow $T_*=0$ and $\tau_0=0$ by letting $\eta\equiv1$. In this case, \mref{p0} now is
$$\barr{ll}\lefteqn{\sup_{t\in [0,T']}\iidx{\Og}{|u|^{l+2}} + \llg_l\itQ{Q}{\llg(u) |u|^l|Du|^2} \le}\hspace{5cm}&\\& C\itQ{Q}{\myprod{\hat{f}(u,Du),|u|^lu}}+\iidx{\Og}{|u(x,0)|^{l+2}}.\earr$$ We can see that the proof can continue to give \mref{p00z} with the right hand side depending on $\|u(\cdot,0)\|_{L^{l+2}(\Og)}$. This suffices to give the global existence result of  \reftheo{thmglobext}. Once this is established, we can take $\tau_0=1$ to get a uniform bound for $\|u(\cdot,t)\|_{L^{l+2}(\Og)}$, independent of $\|u(\cdot,0)\|_{L^{l+2}(\Og)}$, and \reftheo{thmglobatt} can apply.
\erem

\brem{fduremz} By the ellipticity of $A(u)$ \mref{sgineq} holds if $l=0$, i.e. we test the system with $u$, and we do not need SG') here. In this case, \mref{llgDu2} provides the estimate  $$ \int_{T+\tau_0}^{T'}\iidx{\Og}{\llg(u)|Du|^2\eta^p}dxdt \le C(\sup_{(\tau_0,T')}C_0(t),\tau_0)|\Og|(T'-T).$$ This is \mref{llgDu2z} which was used in the proof of \reflemm{ldulemm1z}.
\erem

We conclude this paper by giving a simple example proving the existence of a global attrator of the generalized version \mref{e0}.
Inspired by the {\em competitive} Lotka-Volterra reaction in \mref{e0}, we assume that $\hat{f}(u,Du)$ is of the form
\beqno{sktf}\hat{f}(u,Du)=B(u)Du+Ku -G(u)u,\eeq where $K,B(u),G(u)$ are $m\times m$ matrices and $K$ is a constant one. 

\btheo{SKTyag} Assume SKT), \mref{sktf} and that
$B(u),G(u)$ are $C^1$ in $u$. We assume that $|B(u)|\le C\llg^\frac12(u)$ for all $u\in\RR^m$. In addition, $G(u)$ is  positive definite in the following sense: there are $c_0>0$ and $\kappa\in(0,k]$ such that for all $w,u\in\RR^m$ \beqno{Gcond}\myprod{G(w)u,u} \ge  c_0|w|^{\kappa}|u|^2, \; |G(u)|\sim |u|^\kappa, \; |G_u(u)|\sim |u|^{\kappa-1}.\eeq

Then \mref{e00} defines a dynamical system which possesses a global attractor on $X=W^{1,p_0}(\Og)$, $p_0>2$.
\etheo

\bproof We first consider the global existence by applying \reftheo{thmglobext} (and then \reftheo{L1thm}).  By \refrem{fsigmarem}, we can replace $\hat{f}(\sg u, \sg Du)$ in \mref{e1famzz} by $$\hat{f}_{\sg}(u,Du)= \sg^kB(\sg u)Du+f_{\sg}(u),\; f_{\sg}(u):=\sg^k Ku - \sg^{k-\kappa}G(\sg u)u\quad \sg\in[0,1].$$

Since  $|\partial _u f_\sg(u)|\le \sg^k|K|+\sg^{k-\kappa}|G(\sg u)|+ \sg^{k-\kappa+1}|u||\partial_{\sg u}G(\sg u)|\le C\llg(\sg u)$ by the last two conditions in \mref{Gcond} and the fact that $\llg(\sg u)\sim \llg_0+\sg^k|u|^k$, we see that $f_\sg$ satisfy \mref{fllg} in F'). Hence, the argument in the proof of lemmas leading to \reftheo{thmglobext} continues to hold if we can show that $\|u_\sg\|_{L^1(\Og)}$ can be bounded uniformly for all solutions $u_\sg$ to \mref{e1famzz} in any finite time intervals. 

For any $\sg\in[0,1]$, let $u$ be the solution of \mref{e1famzz}. From \mref{Gcond}, with $w=\sg u$, the bound on $B(u)$ and Young's inequality, we have
$$\myprod{G(\sg u)u,u} \ge C_1\sg^\kappa|u|^{\kappa+2},\; \myprod{B(\sg u)Du,u}\le \eg\llg(\sg u)|Du|^2 + C(\eg)|u|^2.$$

Multiplying \mref{e1famzz} with $u$, integrating by parts in $x$, and using the above inequalities with sufficiently small $\eg$, we easily obtain $$\frac{d}{dt}\iidx{\Og}{|u|^2}+\iidx{\Og}{\llg(\sg u)|Du|^2}\le C_1\sg^{k}\iidx{\Og}{|u|^2}-C_2\sg^{k}\iidx{\Og}{|u|^{\kappa+2}}.$$ 
By H\"older's inequality, we can find a constant $C_3>0$ such that
\beqno{holderp}C_3\left(\iidx{\Og}{|u|^2}\right)^{p}\le C_2\iidx{\Og}{|u|^{\kappa+2}}, \quad p:=(\kappa+2)/2>1.\eeq

Therefore, for $y(t):= \|u(\cdot,t)\|_{L^2(\Og)}^2$ $$y'(t)\le F(y(t)), \quad F(y):=\sg^k(C_1y-C_3y^p).$$

Since $p>1$, we see that $F(y)\le0$ if $y\ge y_* :=(C_1/C_3)^{1/(p-1)}$. As $y(0)\ge 0$, it follows that $y(t)\le \max\{y(0),y_*\}$. Hence, $\|u\|_{L^2(\Og)}^2$ is bounded by a constant, independently of $\sg$. This gives a uniform bound for solutions to \mref{e1famzz}. The global existence result then follows.

We now turn to the existence of global attractors. We test the system with $u$ and obtain, denoting $|K|$ the matrix norm of $K$
\beqno{K}\frac{d}{dt}\iidx{\Og}{|u|^2}
+\iidx{\Og}{\llg(u)|Du|^2}\le \iidx{\Og}{(C\llg^\frac12(u)|Du||u|+|K||u|^2-c_0|u|^{\kappa+2})}.\eeq

Using Young's inequality, for any $\eg>0$ we can find a constant $C(\eg)$ such that $$C\llg^\frac12(u)|Du||u|\le \eg\llg(u)|Du|^2+C(\eg)|u|^2,\;|u|^2\le \eg|u|^{\kappa+2}+ C(\eg).$$ Applying the above inequalities to the first and second integrands on the right hand side of \mref{K}, for sufficiently small $\eg$ we then deduce  the following.
$$\frac{d}{dt}\iidx{\Og}{|u|^2}+c_1 \iidx{\Og}{|u|^{\kappa+2}}\le c_2.$$
Here, $c_1,c_2$ are postive constants depending only on $|K|,c_0$.  As in \mref{holderp}, we can apply H\"older's inequality to the second term on the left hand side to get another positive constant $c_3$ which depend only on $|K|,c_1,|\Og|$ such that for $y(t):=\|u(\cdot,t)\|_{L^2(\Og)}$ and all $t>0$ 
$$y'+c_3y^{p}\le c_2, \quad p=(\kappa+2)/2>1.$$

Using the uniform Gronwall's lemma (\cite[Lemma 5.1]{tem}) with $\cg=c_3$ and $\dg=c_2$, we have \beqno{ybd1}y(t)\le (c_2/c_3)^{1/p} + (c_3(p-1)t)^{-1/(p-1)}.\eeq

For any fixed $M_1>(c_2/c_3)^{1/p}$ we let $$T_*=\frac{1}{c_3(p-1)}\left(M_1-(c_2/c_3)^{1/p}\right)^{1-p}.$$ It is easy to see from \mref{ybd1} that $y(t)\le M_1$ if $t\ge T_*$. The existence of the golbal attractor follows from this uniform estimate. \eproof

We remark that \reftheo{SKTyag} applies to the system \mref{e0} with competitive Lotka-Volterra reaction terms and positive initial data $U_0$. However, the parameters $A(u),f(u)$ of this system satisfy the assumption A),F) only for positive $u$. The positivity of solutions to \mref{e0} was established in \cite{yag} under suitable conditions on the parameter $\ag_{ij}$'s. In a forthcoming paper \cite{LDpos} we will show that this is the case even in a much more general setting.

\bibliographystyle{plain}

\end{document}